\documentclass[11pt]{article}
\usepackage{amssymb,amsmath}
\usepackage{graphicx} %%%%!!!
\newcommand{\cl}{\mathop{\rm cl}}
\newcommand{\Proj}{\mathop{\rm Proj_{rot}}}

\newcommand{\inter}{\mathop{\rm Int}}

\numberwithin{equation}{section}
\newtheorem{theorem}{Theorem}[section]
\newtheorem{corollary}{Corollary}[section]

\title{Rigidity Conditions for the Boundaries of Submanifolds in a
Riemannian Manifold\footnote{
Riemannian manifold, intrinsic metric, induced boundary metric,
strict convexity of submanifold, geodesics, rigidity conditions}}
\author{Anatoly~P.~Kopylov\footnote{Sobolev Institute of Mathematics,
Acad. Koptyuga pr. 4, and Novosibirsk State University, Pirogova
str., 2, 630090 Novosibirsk, Russia; apkopylov@yahoo.com},
Mikhail~V.~Korobkov\footnote{Sobolev Institute of Mathematics,
Acad. Koptyuga pr. 4, and Novosibirsk State University, Pirogova
str., 2, 630090 Novosibirsk, Russia; korob@math.nsc.ru}}
\date{01.01.2014}

\begin{document}

\maketitle

%\section{Abstract}\label{s0}
\begin{abstract}

Developing A.~D.~Aleksandrov's ideas, in~\cite{Ko}
(see also~\cite{Ko1}), the first-named author of this article proposed
the following approach to study of rigidity problems for the
boundary of a
$C^0$-submanifold
in a smooth Riemannian manifold: Let
$Y_1$
be a
$2$-dimensional
compact connected
$C^0$-submanifold
with nonempty boundary in a
$2$-dimensional
smooth connected Riemannian manifold
$(X,g)$
without boundary satisfying the condition
$\rho_{Y_1}(x,y) = \liminf_{x' \to x, y' \to y, x',y' \in \inter Y_1}
\{[l(\gamma_{x',y',\inter Y_1})]\} < \infty$,
if
$x,y \in Y_1$.
Here
$\inf[l(\gamma_{x',y',\inter Y_1})]$
is the infimum of the length of smooth paths joining
$x'$
and
$y'$
in the interior
$\inter Y_1$
of
$Y_1$.
In the present paper, we first establish that
$\rho_{Y_1}$
is a metric on
$Y_1$.
Suppose further that
$Y_1$
is strictly convex in the metric
$\rho_{Y_1}$
(see Sec.~\ref{s3}). Consider another
$2$-dimensional
compact connected
$C^0$-submanifold
$Y_2$
of
$X$
with boundary satisfying the condition
$\rho_{Y_2}(x,y) < \infty$,
$x,y \in Y_2$,
and assume that
$\partial Y_1$
and
$\partial Y_2$
are isometric in the metrics
$\rho_{Y_j}$,
$j = 1,2$.
There appears the following natural question: Under which additional
conditions are the boundaries
$\partial Y_1$
and
$\partial Y_2$
of
$Y_1$
and
$Y_2$
isometric in the metric
$\rho_X$
of the ambient manifold
$X$?
The paper is devoted to the detailed discussions of this question.
In it, we in particular obtain new results
concerning the rigidity problems for the boundaries of
$C^0$-submanifolds
in a Riemannian manifold. The case of
$\dim Y_j = \dim X = n$,
$n > 2$,
is also considered.

\end{abstract}

%\medskip

%\tableofcontents

\section{Introduction: Unique Determination of Surfaces by Their
Relative Metrics on Boundaries}\label{s1}

A classical theorem says (see~\cite{Po}): \textit{
If two bounded closed convex surfaces in the three-dimensional
Euclidean space are isometric in their intrinsic
metrics  then they are equal{\rm,} i.e.{\rm,} they can be matched
by a motion.}

The problems of unique determination of closed convex surfaces by
their intrinsic metrics goes back to the result of Cauchy,
obtained already in 1813, that \textit{any closed convex
polyhedrons
$P_1$
and
$P_2$
{\rm(}in the three-dimensional Euclidean space{\rm)} that are
equally composed of congruent faces are equal.} Since then this
problem has been studied by many people for about 140 years (for example,
by Minkowski, Hilbert, Weyl, Blaschke, Cohn-Vossen, Aleksandrov,
Pogorelov and other prominent mathematicians (see, for instance,
the historical overview in~\cite{Po}, Chapter 3); finally, its
complete solution, which is just the theorem we have cited at
the beginning, was obtained by A.~V.~Pogorelov. For generalizations
of Pogorelov's result to higher dimensions, see~\cite{Se}.

In~\cite{Ko2}, we proposed a new approach to the problem of
unique determination of surfaces, which enabled us to substantially
enlarge the framework of the problem. The following model
situation illustrates the essence of this approach fairly well:

Let
$U_1$
and
$U_2$
be two domains (i.e., open connected sets) in the real
$n$-dimensional
Euclidean space
$\mathbb R^n$
whose closures
$\cl U_j$,
where
$j = 1,2$,
are Lipschitz manifolds (such that
$\partial (\cl U_j) = \partial U_j \ne\varnothing$,
where
$\partial E$
is the boundary of
$E$
in
$\mathbb R^n$).
Assume also that the boundaries
$\partial U_1$
and
$\partial U_2$
of these domains, which coincide with the boundaries of the
manifolds
$\cl U_1$
and
$\cl U_2$,
are isometric with respect to their relative metrics
$\rho_{\partial U_j,U_j}$
($j = 1,2$),
i.e., with respect to the metrics that are the restrictions to
the boundaries
$\partial U_j$
of the extensions
$\rho_{\cl U_j}$
(by continuity) of the intrinsic metrics
$\rho_{U_j}$
of the
domains
$U_j$
to
$\cl U_j$.
The following natural question arises: \textit{Under which
additional conditions are the domains
$U_1$
and
$U_2$
themselves isometric {\rm(}in the Euclidean metric{\rm)}}?
In particular, the natural character of this problem is
determined by the circumstance that the problem of unique
determination of closed convex surfaces mentioned at the
beginning of the article is its most important particular case.
Indeed, assume that
$S_1$
and
$S_2$
are two closed convex surfaces in
$\mathbb R^3$,
i.e., they are the boundaries of two bounded convex domains
$G_1 \subset \mathbb R^3$
and
$G_2 \subset \mathbb R^3$.
Let
$U_j = \mathbb R^3 \setminus \cl G_j$
be the complement of the closure
$\cl G_j$
of the domain
$G_j$,
$j = 1,2$.
Then the intrinsic metrics on the surfaces
$S_1 = \partial U_1$
and
$S_2 = \partial U_2$
coincide with the relative metrics
$\rho_{\partial U_1,U_1}$
and
$\rho_{\partial U_2,U_2}$
on the boundaries of the domains
$U_1$
and
$U_2$,
and thus the problem of unique determination of closed convex
surfaces by their intrinsic metrics is indeed a particular case
of the problem of unique determination of domains by the
relative metrics on their boundaries.

The generalization of the problem of the unique determination of
surfaces ensuing from a new approach suggested in~\cite{Ko2}
manifests itself in the fact that the unique determination of
domains by the relative metrics on their boundaries holds not
only when their complements are bounded convex sets but, for
example, also in the following cases.

\textit{The domain
$U_1$
is bounded and convex and the domain
$U_2$
is arbitrary} (A.~P.~Kopylov (see~\cite{Ko2})).

\textit{The domain
$U_1$
is strictly convex and the domain
$U_2$
is arbitrary} (A.~D.~Aleksandrov (see~\cite{VA1})).

\textit{The domains
$U_1$,
$U_2$
are bounded and their boundaries are smooth} (V.~A.~Aleksandrov
(see~\cite{VA1})).

\textit{The domains
$U_1$
and
$U_2$
have nonempty bounded complements, while their boundaries are
($n - 1$)-dimensional
connected
$C^1$-manifolds
without boundary{\rm,}
$n > 2$} (V.~A.~Aleksandrov (see~\cite{VA2})).

In the papers~\cite{Kor1}-\cite{Kor3}, M.~V.~Korobkov
(in particular) obtained a complete solution to the problem of unique
determination of a plane (space) domain in the class of all
plane (space) domains by the relative metric on its boundary.

In this connection, there appears the following question: Is it possible
to construct an analog of the theory of rigidity of surfaces in
Euclidean spaces in the general case of the boundaries of
submanifolds in Riemannian manifolds?

Our article is devoted to a detailed discussion of this question.
In it, we in particular obtain
new results concerning rigidity problems for the boundaries of
$n$-dimensional
connected submanifolds with boundary in
$n$-dimensional
smooth connected Riemannian manifolds without boundary
($n \ge 2$).

In what follows, all paths
$\gamma: [\alpha,\beta] \to \mathbb R^n$,
where
$\alpha,\beta \in \mathbb R$,
are assumed continuous and non-constant, and
$l(\gamma)$
means the length of a path
$\gamma$.

\section{Rigidity Problems and Intrinsic Geometry of Submanifolds
in Riemannian Manifolds}\label{s2}

Let
$(X,g)$
be an
$n$-dimensional
smooth connected Riemannian manifold without boundary and let
$Y$
be an
$n$-dimensional
compact connected
$C^0$-submanifold
in
$X$
with nonempty boundary
$\partial Y$
($n \ge 2$).

A classical object of investigations (see, for example,
~\cite{Al}) is given by the intrinsic metric
$\rho_{\partial Y}$
on the hypersurface
$\partial Y$
defined for
$x,y \in \partial Y$
as the infimum of the lengths of curves
$\nu \subset \partial Y$ joining
$x$
and
$y$.
In the recent decades, an alternative approach arose in the
rigidity theory for submanifolds of Riemannian manifolds (see,
for instance, the recent articles~\cite{Kor3},~\cite{Ko},
and~\cite{Ko1}, which also contain a historical survey of works
on the topic). In accordance with this approach, the metric on
$\partial Y$
is induced by the intrinsic metric of the interior
$\inter Y$
of the submanifold
$Y$.

Namely, suppose that
$Y$
satisfies the following condition:

${\rm(i)}$
if
$x,y \in Y$,
then
\begin{equation}\label{eq2.1}
\rho_Y(x,y) = \liminf_{x' \to x, y' \to y; x',y' \in \inter Y}
\{\inf[l(\gamma_{x',y',\inter Y})]\} < \infty,
\end{equation}
where
$\inf[l(\gamma_{x',y',\inter Y})]$
is the infimum of the lengths
$l(\gamma_{x',y',\inter Y})$
of smooth paths
$\gamma_{x',y',\inter Y} : [0,1] \to \inter Y$
joining
$x'$
and
$y'$
in the interior
$\inter Y$
of Y.

\textbf{Remark~2.1.}
Easy examples show that if
$X$
is an
$n$-dimensional
connected smooth Riemannian manifold without boundary then an
$n$-dimensional
compact connected
$C^0$-submanifold
in
$X$
with nonempty boundary may fail to satisfy condition
${\rm(i)}$.
For
$n = 2$,
we have the following counterexample:

Let
$(X,g)$
be the space
$\mathbb R^2$
endowed with the Euclidean metric and let
$Y$
be a closed Jordan domain in
$\mathbb R^2$
whose boundary is the union of the singleton
$\{0\}$
consisting of the origin
$0$,
the segment
$\{(1 - t)(e_1 + 2e_2) + t(e_1 + e_2) : 0 \le t \le 1\}$,
and
of the segments of the following four types:
$$
\biggl\{\frac{(1 - t)(e_1 + e_2)}{n} + \frac{te_1}{n + 1} :
0 \le t \le 1\biggr\}
\quad (n = 1,2,\dots);
$$
$$
\biggl\{\frac{e_1 + (1 - t)e_2}{n} : 0 \le t \le 1\biggr\}
\quad (n = 2,3,\dots);
$$
$$
\biggl\{\frac{(1 - t)(e_1 + 2e_2)}{n} +
\frac{2t(2e_1 + e_2)}{4n + 3} : 0 \le t \le 1\biggr\} \quad
(n = 1,2,\dots);
$$
$$
\biggl\{\frac{(1 - t)(e_1 + 2e_2)}{n + 1} +
\frac{2t(2e_1 + e_2)}{4n + 3} : 0 \le t \le 1\biggr\} \quad
(n = 1,2,\dots).
$$
Here
$e_1,e_2$ is the canonical basis in
$\mathbb R^2$.
By the construction of
$Y$,
we have
$\rho_Y(0,E) = \infty$
for every
$E \in Y \setminus \{0\}$
(see figure~1).

\begin{center}     %%%%!
                   %%%%!
\begin{figure}[h]  %%%%!
\includegraphics{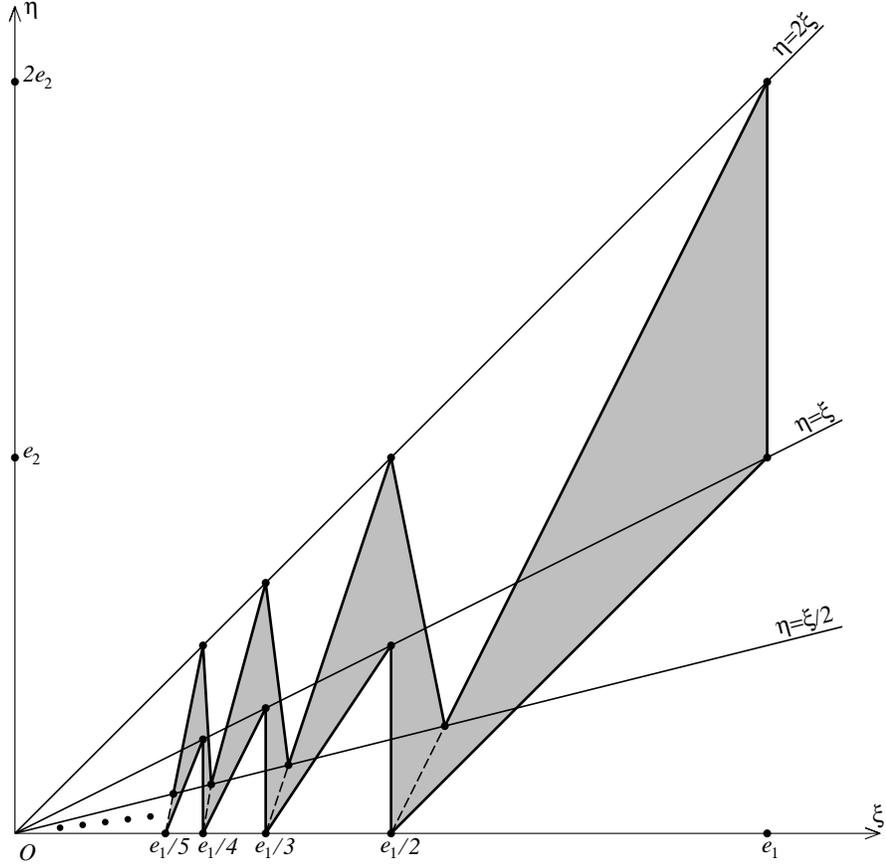}  %%%%!
\caption{An example of
$2$-dimensional
compact connected
$C^0$-submanifold
with nonempty boundary which does not satisfy condition
${\rm(i)}$}  %%%%!
\end{figure}       %%%%!
                   %%%%!
\end{center}       %%%%!

\textbf{Remark~2.2.}
Note that if
$X = \mathbb R^n$
and
$U$
is a domain in
$\mathbb R^n$
whose closure
$Y = \cl U$
is a Lipschitz manifold (such that
$\partial(\cl U) = \partial U \ne\varnothing$),
then
$\rho_{\partial U,U}(x,y) = \rho_Y(x,y)$
($x,y \in \partial U$)
and
$Y$
satisfies~
${\rm(i)}$.
Hence, this example is an important particular case of
submanifolds
$Y$
in a Riemannian manifold
$X$
satisfying
${\rm(i)}$.

To prove our rigidity results for boundaries of submanifolds in
a Riemannian manifold (see Sec.~\ref{s3}), we first need to study the
properties of the intrinsic geometry of these submanifolds.

One of the main results of this section is as follows:

\begin{theorem}\label{t2.1}
Let
$n = 2$.
Then, under condition
${\rm(i)}$,
the function
$\rho_Y$
defined by {\rm(\ref{eq2.1})}
is a metric on
$Y$.
\end{theorem}

\textit{Proof.} It suffices to prove that $\rho_Y$ satisfies the
triangle inequality. Let $A$, $O$, and $D$ be three points on the
boundary of $Y$ (note that this case is basic because the other
cases are simpler). Consider $\varepsilon > 0$ and assume that
$\gamma_{A_{\varepsilon} O^1_{\varepsilon}} : [0,1] \to \inter Y$
and $\gamma_{O^2_{\varepsilon} D_{\varepsilon}} : [2,3] \to \inter
Y$ are smooth paths with the endpoints $A_{\varepsilon} \gamma_{A_{\varepsilon} O^1_{\varepsilon}}(0)$, $O^1_{\varepsilon}
= \gamma_{A_{\varepsilon} O^1_{\varepsilon}}(1)$ and
$D_{\varepsilon} = \gamma_{O^2_{\varepsilon} D_{\varepsilon}}(3)$,
$O^2_{\varepsilon} = \gamma_{O^2_{\varepsilon}
D_{\varepsilon}}(2)$ satisfying the conditions
$\rho_X(A_{\varepsilon}, A) \le \varepsilon$,
$\rho_X(D_{\varepsilon}, D) \le \varepsilon$,
$\rho_X(O^j_{\varepsilon}, O) \le \varepsilon$ ($j = 1;\; 2$),
$|l(\gamma_{A_{\varepsilon} O^1_{\varepsilon}}) - \rho_Y(A,O)| \le
\varepsilon$, and $|l(\gamma_{O^2_{\varepsilon} D_{\varepsilon}})
- \rho_Y(O,D)| \le \varepsilon$. Let $(U, h)$ be a chart of the
manifold $X$ such that $U$ is an open neighborhood of the point
$O$ in $X$, $h(U)$ is the unit disk $B(0, 1) = \{(x_1, x_2) \in
\mathbb R^2 : x_1^2 + x_2^2 < 1\}$ in $\mathbb R^2$, and $h(O) 0$ ($0 = (0, 0)$ is the origin in $\mathbb R^2$); moreover $h : U
\to h(U)$ is a diffeomorphism having the following property: there
exists a chart $(Z,\psi)$ of $Y$ with $\psi(O) = 0$, $A,D \in U
\setminus \cl_X Z$ ($\cl_X Z$ is the closure of $Z$ in the space
$(X,g)$) and $Z = \widetilde{U} \cap Y$ is the intersection of an
open neighborhood $\widetilde{U}$ ($\subset U$) of $O$ in $X$ and
$Y$ whose image $\psi(Z)$ under $\psi$ is the half-disk $B_+(0, 1)
= \{(x_1, x_2) \in B(0, 1) : x_1 \ge 0\}$. Suppose that $\sigma_r$
is an arc of the circle $\partial B(0, r)$ which is a connected
component of the set $V \cap \partial B(0, r)$, where $V = h(Z)$
and $0 < r < r^* = \min \{|h(\psi^{-1}(x_1, x_2))| : x_1^2 + x_2^2
= 1/4, \, x_1 \ge 0\}$. Among these components, there is at least
one (preserve the notation $\sigma_r$ for it) whose ends belong to the
sets $h(\psi^{-1}(\{-te_2 : 0 < t < 1\}))$ and $h(\psi^{-1}(\{te_2
: 0 < t < 1\}))$ respectively. Otherwise, the closure of the
connected component of the set $V \cap B(0, r)$ whose boundary
contains the origin would contain a point belonging to the arc
$\{e^{i\theta}/2 : |\theta| \le \pi/2\}$ (here we make use of the
complex notation $z = re^{i\theta}$ for points $z \in \mathbb R^2$
($= \mathbb C$)). But this is impossible. Therefore, the
above-mentioned arc $\sigma_r$ exists.

It is easy to check that if $\varepsilon$ is sufficiently small
then the images of the paths $h \circ
\gamma_{A_{\varepsilon}O^1_{\varepsilon}}$ and $h \circ
\gamma_{O^2_{\varepsilon}D_{\varepsilon}}$ also intersect the arc
$\sigma_r$, i.e., there are $t_1 \in ]0,1[$, $t_2 \in ]2,3[$ such
that $\gamma_{A_\varepsilon O^1_\varepsilon}(t_1) = x^1 \in Z$,
$\gamma_{O^2_\varepsilon D_\varepsilon}(t_2) = x^2 \in Z$ and
$h(x^j) \in \sigma_r$, $j = 1,2$. Let $\gamma_r : [t_1,t_2] \to
\sigma_r$ be a smooth parametrization of the corresponding subarc
of $\sigma_r$, i.e., $\gamma_r(t_j) = h(x^j)$, $j = 1,2$. Now we
can define a mapping $\widetilde{\gamma}_\varepsilon : [0,3] \to \inter Y$ by
setting
$$\widetilde{\gamma}_\varepsilon(t)=\left\{ \aligned
\gamma_{A_{\varepsilon} O^1_{\varepsilon}}(t) , & \ \ t\in [0,t_1]; \\
h^{-1}(\gamma_r(t)), & \ \ t\in]t_1,t_2[; \\
\gamma_{O^2_{\varepsilon}D_{\varepsilon} }(t) , & \ \ t\in
[t_2,3].
\endaligned \right.
$$
By construction, $\widetilde{\gamma}_\varepsilon$ is a piecewise smooth path
joining the points $A_\varepsilon = \widetilde{\gamma}_\varepsilon(0)$,
$D_\varepsilon = \widetilde{\gamma}_\varepsilon(3)$ in $\inter Y$; moreover,
$$
l(\widetilde{\gamma}_\varepsilon) \le l(\gamma_{A_\varepsilon
O^1_\varepsilon}) + l(\gamma_{O^2_\varepsilon D_\varepsilon}) +
l(h^{-1}(\sigma_r)).
$$
By an appropriate choice of $\varepsilon > 0$, we can make~$r
> 0$ arbitrarily small, and since a piecewise smooth
path can be approximated by smooth paths, we have $\rho_Y(A,D) \le
\rho_Y(A,O) + \rho_Y(O,D)$, q.e.d.

\bigskip

In connection with Theorem~\ref{t2.1}, there appears a natural question:
Are there analogs of this theorem for
$n \ge 3$?
The following Theorem~\ref{t2.2} answers this question in the negative:

\begin{theorem}\label{t2.2}
If
$n \ge 3$
then there exists an
$n$-dimensional
compact connected
$C^0$-manifold
$Y \subset \mathbb R^n$
with nonempty boundary
$\partial Y$
such that condition
${\rm(i)}$
{\rm(}where now
$X = \mathbb R^n)$
is fulfilled for
$Y$
but the function
$\rho_Y$
in this condition is not a metric on
$Y$.
\end{theorem}

\textit{Proof.} It suffices to consider the case of $n = 3$. Suppose that $A$, $O$,
$D$ are points in $\mathbb R^3$, $O$ is the origin in $\mathbb
R^3$, $|A| = |D| = 1$, and the angle between the segments $OA$ and
$OD$ is equal to $\frac{\pi}{6}$.

The manifold $Y$ will be constructed so that $O \in \partial Y$,
and $]O,A] \subset \inter Y$, $]O,D] \subset \inter Y$. Under
these conditions, $\rho_Y(O,A) = \rho_Y(O,D) =1$. However, the
boundary of $Y$ will create ``obstacles'' between $A$ and $D$ such
that the length of any curve joining $A$ and $D$ in $\inter Y$
will be greater than $\frac{12}{5}$ (this means the violation of
the triangle inequality for~$\rho_Y$).

Consider a countable collection of mutually disjoint segments
$\{I^k_j\}_{j \in \mathbb N,\, k = 1,\dots,k_j}$ lying in the
interior of the triangle $6 \Delta AOD$ (which is obtained from
the original triangle $\Delta AOD$ by dilation with
coefficient~$6$) with the following properties:

$(*)$ every segment $I^k_j =[x^k_j,y^k_j]$ lies on a ray starting
at the origin, $y^k_j = 11 x^k_j$, and $|x^k_j| = 2^{-j}$;

$(**)$ any curve $\gamma$ with ends $A$, $D$ whose interior
points lie in the interior of the triangle $4 \Delta AOD$ and
belong to no~segment $I^k_j$, satisfies the estimate $l(\gamma) \ge 6$.

The existence of such a family of segments is certain:
the segments of the family must be situated chequerwise
so that any curve disjoint from them be
sawtooth, with the total length of its ``teeth'' greater than $6$
(it can clearly be made greater than any prescribed positive
number). However, below we exactly describe the construction.

It is easy to include the above-indicated family of segments in
the boundary $\partial Y$ of~$Y$. Thus, it creates a desired
``obstacle'' to joining $A$ and $D$ in the plane of $\Delta AOD$.
But it makes no obstacle to joining $A$ and $D$ in the space. The
simplest way to create such a~space obstacle is as follows: Rotate
each segment $I^k_j$ along a spiral around the axis $OA$. Make the
number of coils so large that the length of this spiral be large
and its pitch (i.e., the distance between the origin and the end
of a coil) be sufficiently small. Then the set $S^k_j$ obtained as
the result of the rotation of the segment $I^k_j$ is diffeomorphic
to a plane rectangle, and it lies in a small neighborhood of the
cone of revolution with axis $AO$ containing the segment $I^k_j$.
The last circumstance guarantees that the sets $S^k_j$ are
disjoint as before, and so (as above) it is easy to include them
in the boundary $\partial Y$ but, due to the properties of the
$I^k_j$'s and a~large number of coils of the spirals $S^k_j$, any
curve joining $A,D$ and disjoint from each $S^k_j$ has length $\ge
\frac{12}{5}$.

We turn to an exact description of the constructions used. First
describe the construction of the family of segments $I^k_j$. They
are chosen on the basis of the following observation:

Let $\gamma : [0,1] \to 4 \Delta AOD$ be any curve with ends
$\gamma (0) = A$, $\gamma (1) = D$ whose interior points lie in
the interior of the triangle $4 \Delta AOD$. For $j \in \mathbb
N$, put $R_j = \{x \in 4 \Delta AOD : |x| \in [8 \cdot 2^{-j}, 4
\cdot 2^{-j}]\}$. It is clear that
$$
4 \Delta AOD \setminus \{O\} = \cup_{j \in \mathbb N} R_j.
$$
Introduce the polar system of coordinates on the plane of the
triangle $\Delta AOD$ with center $O$ such that the coordinates of
the points $A,D$ are $r = 1$, $\varphi = 0$ and $r = 1$, $\varphi
= \frac{\pi}{6}$, respectively. Given a point $x \in 6 \Delta
AOD$, let $\varphi_x$ be the angular coordinate of $x$ in
$[0,\frac{\pi}{6}]$. Let $\Phi_j = \{\varphi_{\gamma(t)} :
\gamma(t) \in R_j\}$. Obviously, there is $j_0 \in \mathbb N$ such
that
\begin{equation}
\mathcal H^1(\Phi_{j_0}) \ge 2^{-j_0} \frac{\pi}{6}, \label{2}
\end{equation}
where $\mathcal H^1$ is the Hausdorff $1$-measure. This means
that, while in the layer $R_{j_0}$, the curve $\gamma$ covers the
angular distance $\ge 2^{-j_0} \frac{\pi}{6}$. The segments
$I^k_j$ must be chosen such that~(\ref{2}) together with the
condition
$$
\gamma(t) \cap I^k_j = \varnothing \quad \forall t \in [0,1]\,\,
\forall j \in \mathbb N\,\, \forall k \in \{1,\dots,k_j\}
$$
give the desired estimate $l(\gamma) \ge 6$. To this end, it
suffices to take $k_j = [(2 \pi)^j]$ (the integral part of $(2
\pi)^j$) and
$$
I^k_j = \{x \in 6 \Delta AOD : \varphi_x = k (2 \pi)^{-j}
\frac{\pi}{6},\, |x| \in [11 \cdot 2^{-j},2^{-j}]\},
$$
$k = 1,\dots,k_j$. Indeed, under this choice of the $I^k_j$'s,
estimate~(\ref{2}) implies that $\gamma$ must intersect at least
$(2 \pi)^{j_0} 2^{-j_0} = \pi^{j_0} > 3^{j_0}$ of the figures
$$
U_k = \{x \in R_{j_0} : \varphi_x \in (k(2 \pi)^{-j_0}
\frac{\pi}{6}, (k + 1) (2 \pi)^{-j_0} \frac{\pi}{6})\}.
$$
Since these figures are separated by the segments $I^k_{j_0}$ in
the layer $R_{j_0}$, the curve $\gamma$ must be disjoint from them
each time in passing from one figure to another. The number of
these passages must be at least $3^{j_0} - 1$, and a~fragment of
$\gamma$ of length at least $2 \cdot 3 \cdot 2^{-j_0}$ is required
for each passage (because the ends of the segments $I^k_{j_0}$ go
beyond the boundary of the layer $R_{j_0}$ containing the figures
$U_k$ at distance $3 \cdot 2^{-j_0}$). Thus, for all these
passages, a section of $\gamma$ is spent of~length at least
$$
6 \cdot 2^{-j_0} (3^{j_0} -1) \ge 6.
$$
Hence, the construction of the segments $I^k_j$ satisfying~($*$)--($**$)
is finished.

Let us now describe the construction of the above-mentioned space
spirals.

For $x \in \mathbb R^3$, denote by $\Pi_x$ the plane that passes
through $x$ and is perpendicular to the segment $OA$. On
$\Pi_{x^k_j}$, consider the polar coordinates
$(\rho, \psi)$ with origin at the point of intersection of
$\Pi_{x^k_j}$ and $[O,A]$ (in this system, the point $x^k_j$ has
coordinates $\rho =\rho^k_j$, $\psi = 0$). Suppose that a point
$x(\psi) \in \Pi_{x^k_j}$ moves along an~Archimedes spiral,
namely, the polar coordinates of the point $x(\psi)$ are
$\rho(\psi) = \rho^k_j - \varepsilon_j \psi$, $\psi
\in [0,2 \pi M_j]$, where $\varepsilon_j$ is a small parameter to
be specified below, and $M_j \in \mathbb N$ is chosen so large
that the length of any curve passing between all coils of the
spiral is at least $10$.

Describe the choice of $M_j$ more exactly. To this end, consider
the points $x(2 \pi)$, $x(2 \pi (M_j - 1))$, $x(2 \pi M_j)$, which
are the ends of the first, penultimate, and last coils of the
spiral respectively (with $x(0) = x^k_j$ taken as the starting
point of the spiral). Then $M_j$ is chosen so large that the
following condition hold:

$(*_1)$ {\it The length of any curve on the plane $\Pi_{x^k_j}$
joining the segments $[x^k_j,x(2 \pi)]$ and $[x(2 \pi (M_j -
1)),x(2 \pi M_j)]$ and disjoint from the spiral $\{x(\psi) : \psi
\in [0,2 \pi M_j]\}$ is at least $10$.}

Figuratively speaking, the constructed spiral bounds
a~``labyrinth'', the mentioned segments are the entrance to and
the exit from this labyrinth, and thus any path through the
labyrinth has length $\ge 10$.

Now, start rotating the entire segment $I^k_j$ in space along the
above-mentioned spiral, i.e., assume that
$I^k_j(\psi) = \{y \lambda x(\psi) : \lambda \in[1,11]\}$. Thus, the segment
$I^k_j(\psi)$ lies on the ray joining $O$ with $x(\psi)$ and has
the same length as the original segment $I^k_j = I^k_j(0)$. Define
the surface $S^k_j = \cup_{\psi \in [0, 2 \pi M_j]} I^k_j(\psi)$.
This surface is diffeomorphic to a plane rectangle (strip). Taking
$\varepsilon_j > 0$ sufficiently small, we may assume without loss
of generality that $2 \pi M_j \varepsilon_j$ is substantially less
than $\rho^k_j$; moreover, that the surfaces $S^k_j$ are mutually
disjoint (obviously, the smallness of $\varepsilon_j$ does not
affect property~$(*_1)$ which in fact depends on $M_j$).

Denote by $y(\psi)=11x(\psi)$ the second end of the segment
$I^k_j(\psi)$. Consider the trapezium $P^k_j$ with vertices
$y^k_j$, $x^k_j$, $x(2 \pi M_j)$, $y(2 \pi M_j)$ and sides
$I^k_j$, $I^k_j(2 \pi M_j)$, $[x^k_j,x(2 \pi M_j)]$, and
$[y^k_j,y(2 \pi M_j)]$ (the last two sides are parallel since they
are perpendicular to the segment $AO$). By construction, $P^k_j$
lies on the plane $AOD$; moreover, taking $\varepsilon_j$
sufficiently small, we can obtain the situation where the
trapeziums $P^k_j$ are mutually disjoint (since $P^k_j \to I^k_j$
under fixed $M_j$ and $\varepsilon_j \to 0$). Take an arbitrary
triangle whose vertices lie on $P^k_j$ and such that one of these
vertices is also a vertex at an~acute angle in~$P^k_j$. By
construction, this acute angle is at least $\frac{\pi}{2} - \angle
AOD = \frac{\pi}{3}$. Therefore, the ratio of the side of
the~triangle lying inside the trapezium $P^k_j$ to the sum of the
other two sides (lying on the corresponding sides of $P^k_j$) is
at least $\frac{1}{2} \sin \frac{\pi}{3} > \frac{2}{5}$. If we
consider the same ratio for the case of a~triangle with a vertex
at an~obtuse angle of $P^k_j$ then it is greater than
$\frac{1}{2}$. Thus, we have the following property:

$(*_2)$ {\it For arbitrary triangle whose vertices lie on the
trapezium $P^k_j$ and one of these vertices is also a~vertex in
$P^k_j$, the sum of lengths of the sides situated on the
corresponding sides of $P^k_j$ is less than $\frac{5}{2}$ of the
length of the third side} ({\it lying inside} $P^k_j$).

Let a point $x$ lie inside the cone $K$ formed by the rotation of
the angle $\angle AOD$ around the ray $OA$. Denote by $\Proj x$
the point of the angle $\angle AOD$ which is the image of $x$
under this rotation. Finally, let $K_{4 \Delta AOD}$ stand for the
corresponding truncated cone obtained by the rotation of the
triangle $4 \Delta AOD$, i.e., $K_{4 \Delta AOD} = \{x \in K :
\Proj x \in 4 \Delta AOD\}$.

The key ingredient in the proof of our theorem is the following
assertion:

$(*_3)$ {\it For arbitrary space curve $\gamma$ of length less
than $10$ joining the points $A$ and $D$, contained in the
truncated cone $K_{4 \Delta AOD} \setminus \{O\}$, and disjoint
from each strip $S^k_j$, there exists a plane curve
$\tilde{\gamma}$ contained in the triangle $4 \Delta AOD \setminus
\{O\}$ that joins $A$ and $D$ is disjoint from all segments
$I^k_j$ and such that the length of $\tilde{\gamma}$ is  less than
 $\frac{5}{2}$ of the length of} $\Proj \gamma$.

Prove~$(*_3)$. Suppose that its hypotheses are fulfilled. In
particular, assume that the inclusion $\Proj \gamma \subset 4
\Delta AOD \setminus \{O\}$ holds. We need to modify $\Proj
\gamma$ so that the new curve be contained in the same set but be
disjoint from each of the $I^k_j$'s. The construction splits into
several steps.

{\bf Step~1.} If $\Proj \gamma$ intersects a segment $I^k_j$ then
it necessarily intersects also at least one of the shorter sides
of $P^k_j$.

Recall that, by construction, $P^k_j = \Proj S^k_j$; moreover,
$\gamma$ intersects no spiral strip $S^k_j$. If $\Proj \gamma$
intersected $P^k_j$ without intersecting its shorter sides then
$\gamma$ would pass through all coils of the corresponding spiral.
Then, by~$(*_1)$, the length of the corresponding fragment of
$\gamma$ would be $\ge 10$ in contradiction to our assumptions.
Thus, the assertion of step~1 is proved.

{\bf Step~2.} Denote by $\gamma_{P^k_j}$ the fragment of the plane
curve $\Proj \gamma$ beginning at the first point of its entrance
into the trapezium $P^k_j$ to the point of its exit from $P^k_j$
(i.e., to its last intersection point with $P^k_j$). Then this
fragment $\gamma_{P^k_j}$ can be deformed without changing the
first and the last points so that the corresponding fragment of
the new curve lie entirely on the union of the sides of $P^k_j$;
moreover, its length is less than $\frac{5}{2}$ of the length of
$\gamma_{P^k_j}$.

The assertion of step~2 immediately follows from the assertions of
step~1 and $(*_2)$.

The assertion of step~2 in turn directly implies the desired
assertion $(*_3)$. The proof of $(*_3)$ is finished.

Now, we are ready to pass to the final part of the proof of
Theorem~\ref{t2.2}.

$(*_4)$ {\it The length of any space curve $\gamma \subset \mathbb
R^3 \setminus \{O\}$ joining $A$ and $D$ and disjoint from each
strip $S^k_j$ is at least} $\frac{12}{5}$.

Prove the last assertion. Without loss of generality, we may also
assume that all interior points of $\gamma$ are inside the cone
$K$ (otherwise the initial curve can be modified without any
increase of its length so that assumptions of~$(*_4)$
are still fulfilled and the modified curve lies in~$K$). If
$\gamma$ is not included in the truncated cone $K_{4 \Delta AOD}
\setminus \{O\}$ then $\Proj \gamma$ intersects the segment
$[4A,4D]$; consequently, the length of $\gamma$ is at least $2(4
\sin \angle OAD - 1) = 2(4 \sin \frac{\pi}{3} - 1) = 2 (2 \sqrt 3
- 1) > 4$, and the desired estimate is fulfilled. Similarly, if
the length of $\gamma$ is at least $10$ then the desired estimate
is fulfilled automatically, and there is nothing to prove. Hence,
we may further assume without loss of generality that $\gamma$ is
included in the truncated cone $K_{4 \Delta AOD} \setminus \{O\}$
and its length is less than $10$. Then, by~$(*_3)$, there is a
plane curve $\tilde \gamma$ contained in the triangle $4 \Delta
AOD \setminus \{O\}$, joining the points $A$ and $D$, disjoint
from each segment $I^k_j$, and such that the length of $\tilde
\gamma$ is at most $\frac{5}{2}$ of the length of $\Proj \gamma$.
By property~$(**)$ of the family of segments $I^k_j$, the length
of $\tilde \gamma$ is at least $6$. Consequently, the length of
$\Proj \gamma$ is at least $\frac{12}{5}$, which implies the
desired estimate. Assertion~$(*_4)$ is proved.

The just-proven property $(*_4)$ of the constructed objects
implies Theorem~\ref{t2.2}. Indeed, since the strips $S^k_j$ are
mutually disjoint and, outside every neighborhood of the origin
$O$, there are only finitely many of these strips, it is easy to
construct a $C^0$-manifold $Y \subset \mathbb R^3$ that is
homeomorphic to a closed ball (i.e., $\partial Y$ is homeomorphic
to a two-dimensional sphere) and has the following properties:

(I) $O \in \partial Y$, $[A,O[ \cup [D,O[ \subset \inter Y$;

(II) for every point $y \in (\partial Y) \setminus \{O\}$, there
exists a neighborhood $U(y)$ such that $U(y) \cap \partial Y$ is
$C^1$-diffeomorphic to the plane square $[0,1]^2$;

(III) $S^k_j \subset \partial Y$ for all $j \in \mathbb N,\, k 1,\dots,k_j$.

The construction of $Y$ with properties~(I)--(III) can be carried
out, for example, as follows: As the surface of the zeroth step,
take a~sphere containing $O$ and such that $A$ and $D$ are inside
the sphere. At the $j$th step, a small neighborhood of the point
$O$ of our surface is smoothly deformed so that the modified
surface is still smooth, homeomorphic to a~sphere, and contains
all strips $S^k_j$, $k = 1,\dots,k_j$. Besides, we make sure that,
at each step, the so-obtained surface be disjoint from the
half-intervals $[A,O[$ and $[D,O[$, and, as above, contain all
strips $S^k_i$, $i \le j$, already included therein. Since the
neighborhood we are deforming contracts to the point $O$ as $j \to
\infty$, the so-constructed sequence of surfaces converges (for
example, in the Hausdorff metric) to a limit surface which is the
boundary of a $C^0$-manifold $Y$ with properties~(I)--(III).

Property~(I) guarantees that $\rho_Y(A,O) = \rho_Y(A,D) = 1$ and
$\rho_Y(O,x) \le 1 + \rho_Y(A,x)$ for all $x \in Y$. Property~(II)
implies the estimate $\rho_Y(x,y) < \infty$ for all $x,y \in Y
\setminus \{O\}$, which, granted the previous estimate, yields
$\rho_Y(x,y) < \infty$ for all $x,y \in Y$. However,
property~(III) and the assertion~$(*_4)$ imply that $\rho_Y(A,D)
\ge \frac{12}{5} > 2 = \rho_Y (A,O) + \rho_Y(A,D)$.
Theorem~\ref{t2.2} is proved. q.e.d.

\bigskip

If
$\rho_Y$
is a metric (the dimension
$n$
($\ge 2$)
is arbitrary) then the question of the existence of geodesics is
solved in the following assertion, which implies that
$\rho_Y$
is the \textit{intrinsic metric} (see, for example, \S 6
in~\cite{Al}).

\begin{theorem}\label{t2.3}
Assume that
$\rho_Y$
is a finite function and is a metric on
$Y$.
Then any two points
$x,y \in Y$
can be joined in
$Y$
by a shortest curve
$\gamma : [0,L] \to Y$
in the metric
$\rho_Y;$
i.e.{\rm,}
$\gamma(0) = x,$
$\gamma(L) = y,$
and
\begin{equation}\label{eq2.1''}
\rho_Y(\gamma(s),\gamma(t)) = t - s, \quad \forall s,t \in [0,L],
\quad s < t.
\end{equation}
\end{theorem}

\textit{Proof.} Fix a pair of distinct points $x,y \in Y$ and put
$L = \rho_Y (x,y)$. Now, take a~sequence of paths $\gamma_j :
[0,L] \to \inter Y$ such that $\gamma_j(0) = x_j$, $\gamma_j(L) = y_j$,
$x_j \to x$, $y_j \to y$, and $l(\gamma_j) \to L$ as $j \to
\infty$. Without loss of generality, we may also assume that the
parametrizations of the curves $\gamma_j$ are their natural
parametrizations up to a~factor (tending to~$1$) and the mappings
$\gamma_j$ converge uniformly to a~mapping $\gamma : [0,L] \to Y$
with $\gamma (0) = x$, $\gamma (L) = y$. By these assumptions,
\begin{equation}
\lim_{j \to \infty} l(\gamma_j|_{[s,t]}) = t - s \quad \forall s,t
\in [0,L], \quad s < t. \label{4} \end{equation}

Take an arbitrary pair of numbers $s,t \in [0,L]$, $s < t$. By
construction, we have the convergence $\gamma_j(s) \in \inter Y
\to \gamma(s)$, $\gamma_j(t) \in \inter Y \to \gamma(t)$ as $j \to
\infty$. From here and the definition of the metric
$\rho_Y(\cdot,\cdot)$ it follows that
$$
\rho_Y(\gamma(s),\gamma(t)) \le \lim_{j \to \infty}
l(\gamma_j|_{[s,t]}).
$$
By~(\ref{4}), \begin{equation} \rho_Y(\gamma(s),\gamma(t)) \le t -
s \quad \forall s,t \in [0,L],\,\, s < t. \label{5}
\end{equation} Prove that~(\ref{5}) is indeed an equality. Assume that
$$
\rho_Y(\gamma(s'),\gamma(t')) < t' -s'
$$
for some $s',t' \in [0,L]$, $s' < t'$. Then, applying the triangle
inequality and then~(\ref{5}), we infer
$$
\rho_Y(x,y) \le \rho_Y(x,\gamma(s')) +
\rho_Y(\gamma(s'),\gamma(t')) + \rho_Y(\gamma(t'),y) < s' + (t' -
s') + (L - t') = L,
$$
which contradicts the initial equality $\rho_Y(x,y) = L$. The
so-obtained contradiction completes the proof of
identity~(\ref{eq2.1''}). q.e.d.

\textbf{Remark~2.3}. Identity~(\ref{eq2.1''}) means that the curve
 of Theorem~\ref{t2.3} is a geodesic
in the metric
$\rho_Y$
 , i.e., the length of its fragment between points
$\gamma(s)$,
$\gamma(t)$
calculated in
$\rho_Y$
is equal to
$\rho_Y(\gamma(s),\gamma(t)) = t - s$.
Nevertheless, if
we compute the length of the above-mentioned fragment of the curve in the
initial Riemannian metric then this length need not coincide with
$t - s$;
only
the easily verifiable estimate
$l(\gamma|_{[s,t]}) \le t - s$
holds (see~(\ref{4}). In the general
case, the equality
$l(\gamma|_{[s,t]}) = t - s$
can only be guaranteed if
$n = 2$
(if
$n \ge 3$
then the corresponding counterexample is constructed by analogy with the
counterexample in the proof of Theorem~\ref{t2.2}, see above). In particular, though,
by Theorem~\ref{t2.3}, the metric
$\rho_Y$
is always intrinsic in the sense of the definitions
in \cite[\S 6]{Al}, the space
$(Y,\rho_Y)$
may fail to be \textit{a space with intrinsic metric}
in the sense of [ibid].

\section{Rigidity Theorems for the Boundaries of Submanifolds in
Riemannian Manifolds}\label{s3}

As in Sec.~\ref{s2}, let
$(X,g)$
be an $n$-dimensional
smooth connected Riemannian manifold without boundary and let
$\rho_X$
be its intrinsic metric (i.e., let
$\rho_X(x,y)$
be the infimum of the lengths
$l(\gamma_{x,y,X})$
of smooth paths
$\gamma_{x,y,X} : [0,1] \to X$
joining
points
$x$
and
$y$
in a manifold
$X$).

Assume that
$Y$
is an
$n$-dimensional
compact connected
$C^0$-submanifold
$Y \subset X$
with nonempty boundary
$\partial Y$
satisfying condition
${\rm(i)}$
in Sec.~\ref{s2}, moreover,
$\rho_Y$
is a metric on
$Y$.
Then
$Y$
is called strictly convex in the metric
$\rho_Y$
if, for any
$\alpha,\beta \in Y$,
any shortest path
$\gamma = \gamma_{\alpha,\beta,Y} : [0,1] \to Y$
between
$\alpha$
and
$\beta$
(in the metric
$\rho_Y$)
satisfies
$\gamma(]0,1[) \subset \inter Y$.

\begin{theorem}\label{t3.1}
Let
$n = 2$.
Assume that condition~
${\rm(i)}$
holds for a
$2$-dimensional
compact connected
$C^0$-submanifold
$Y_1$
with nonempty boundary
$\partial Y_1$
of a
$2$-dimensional
smooth connected Riemannian manifold
$X$
without boundary which is strictly convex in the metric
$\rho_{Y_1}$.
Suppose that
$Y_2 \subset X$
is also a
$2$-dimensional
compact connected
$C^0$-submanifold
of
$X$
with
$\partial Y_2 \ne\varnothing$
satisfying
${\rm(i)};$
moreover{\rm,}
$\partial Y_1$
and
$\partial Y_2$
are isometric in the metrics
$\rho_{Y_j},$
for
$j = 1,2$.
Then{\rm,}
$Y_2$
is strictly convex with respect to
$\rho_{Y_2}$.
\end{theorem}

\textit{Proof.}
Suppose that, for points
$x,y \in Y_2$,
there exists a shortest path
$\gamma_{x,y,Y_2} : [0,1] \to Y_2$
in the metric
$\rho_{Y_2}$
joining
$x$
and
$y$
and such that
$\{\gamma_{x,y,Y_2}(]0,1[)\} \cap \partial Y_2 \ne\varnothing$,
i.e.,
$x' = \gamma_{x,y,Y_2}(t') \in \{\gamma_{x,y,Y_2}(]0,1[) \cap
\partial Y_2\}$
for a point
$t' \in ]0,1[$.
By Theorem~\ref{t2.3} and the fact that
$Y_2$
is a
$2$-dimensional
compact connected
$C^0$-submanifold
in
$X$,
for a sufficiently small neighborhood of
$x'$
in
$Y_2$,
we can find points
$x_0,y_0 \in \partial Y_2$
and a shortest path
$\gamma_{x_0,y_0,Y_2} : [0,1] \to Y_2$
between
$x_0$
and
$y_0$
in the same metric satisfying the condition
$x' \in \{\gamma_{x_0,y_0,Y_2}(]0,1[) \cap \partial Y_2\}$.
Further, we will suppose that
$x = x_0$
and
$y = y_0$.

Now, assume that
$f : \partial Y_1 \to \partial Y_2$
is an isometry of
$\partial Y_1$
and
$\partial Y_2$
in the metrics
$\rho_{Y_1}$
and
$\rho_{Y_2}$
of the boundaries
$\partial Y_1$
and
$\partial Y_2$
of the submanifolds
$Y_1$
and
$Y_2$
of
$X$.
From Theorem~\ref{t2.3}, we have
$$
\rho_{Y_2}(x,x') + \rho_{Y_2}(x',y) = l_1 + l_2 = l \rho_{Y_2}(x,y).
$$
Since
$f$
is an isometry,
$$
\rho_{Y_1}(f^{-1}(x),f^{-1}(x')) + \rho_{Y_1}(f^{-1}(x'),f^{-1}(y))
= \rho_{Y_2}(x,x') + \rho_{Y_2}(x',y).
$$
Next, consider shortest paths
$\gamma_{f^{-1}(x),f^{-1}(x'),Y_1} : [0,1/2] \to Y_1$
and
$\gamma_{f^{-1}(x'),f^{-1}(y),Y_1} : [1/2,1] \to Y_1$
in
$\rho_{Y_1}$
between (respectively)
$f^{-1}(x)$
and
$f^{-1}(x')$
and
$f^{-1}(x')$
and
$f^{-1}(y)$,
and then construct a path
$\gamma : [0,1] \to Y_1$
by setting
$\gamma(t) = \gamma_{f^{-1}(x),f^{-1}(x'),Y_1}(t)$
if
$0 \le t < 1/2$
and
$= \gamma_{f^{-1}(x'),f^{-1}(y),Y_1}(t)$
for
$1/2 \le t \le 1$.
Let
$l_{Y_1}(\delta)$
be the length of a path
$\delta : [0,1] \to Y_1$
in the metric $\rho_{Y_1}$.
Since
$\rho_{Y_1}$
is a metric on
$Y_1$,
it is not difficult to show that
$$
l_{Y_1}(\gamma) \le l_{Y_1}(\gamma_{f^{-1}(x),f^{-1}(x'),Y_1}) +
l_{Y_1}(\gamma_{f^{-1}(x'),f^{-1}(y),Y_1}) = l_1 + l_2.
$$
Hence
$\gamma$
is a shortest path in
$\rho_{Y_1}$
joining
$f^{-1}(x)$
and
$f^{-1}(y)$
in
$Y_1$.
This contradicts the strict convexity of
$Y_1$. The theorem is proved.

\begin{corollary}\label{c3.1}
Suppose that the conditions of Theorem~{\rm\ref{t3.1}} hold and
the manifold
$X$
has the following property{\rm:}
$\rho_X(x,y) = \rho_Y(x,y)$
for any two points
$x$
and
$y$
from every
$2$-dimensional
compact connected
$C^0$-submanifold
$Y \subset X$
with
$\partial Y \ne\varnothing$
satisfying condition
${\rm(i)}$
and strictly convex with respect to the metric
$\rho_Y$.
Then,
$\partial Y_1$
and
$\partial Y_2$
are isometric in the metric
$\rho_X$
on the ambient manifold
$X$.
\end{corollary}

\textbf{Remark~3.1.}
The condition imposed on the manifold
$X$
in Corollary~\ref{c3.1} can be reformulated as follows: in this
manifold, every
$2$-dimensional
compact connected
$C^0$-submanifold
$Y$
with boundary satisfying condition
${\rm(i)}$
and strictly convex with respect to its intrinsic metric
$\rho_Y$
is a convex set in the ambient space
$X$
with respect to the metric
$\rho_X$
(for the notion of a convex set in a metric space the reader is referred,
for example, to~\cite{Al}).

We have the following analog of Theorem~\ref{t3.1} and
Corollary~\ref{c3.1} (combined together) for $n \ge 3$:

\begin{theorem}\label{t3.2}
Let
$n \ge 3$.
Suppose that
$(X,g)$
is an
$n$-dimensional
smooth connected Riemannian manifold without boundary and
$Y_1$
and
$Y_2$
are
$n$-dimensional
compact connected
$C^0$-submanifolds with nonempty boundaries
$\partial Y_1$
and
$\partial Y_2$
in
$X$
satisfying conditions
${\rm(i)}$,

${\rm(ii)}$
$\rho_{Y_j}$
is a metric on
$Y_j$
$(j = 1,2)$,

and

${\rm(iii)}$
for any two points
$a,b \in Y_j$,
there exist points
$c,d \in \partial Y_j$
which can be joined in
$Y_j$
by a shortest path
$\gamma : [0,1] \to Y_j$
in the metric
$\rho_{Y_j}$
so that
$a,b \in \gamma([0,1])$.

Furthermore{\rm,} assume that $Y_1$ is strictly convex in the
metric $\rho_{Y_1}$,\, $X$ has the additional property
that $\rho_X(x,y) = \rho_Y(x,y)$ for any two points $x$ and $y$
in every $n$-dimensional compact connected $C^0$-submanifold $Y
\subset X$ with $\partial Y \ne\varnothing$ satisfying
conditions~${\rm(i)}$-${\rm(iii)}$ and strictly convex with respect
to~$\rho_Y$ and the boundaries $\partial Y_1$ and $\partial
Y_2$ of the submanifolds $Y_1$ and $Y_2$ are isometric with
respect to the metrics $\rho_{Y_j},$ where $j = 1,2$. Then,
$\partial Y_1$ and $\partial Y_2$ are isometric with respect to~$\rho_X$.
\end{theorem}

\textbf{Remark~3.2.}
For a submanifold
$Y$
in
$X$,
${\rm(i)}$
and
${\rm(ii)}$
can be considered as conditions of generalized regularity near its boundary.

\textbf{Remark~3.3.}
Theorem~\ref{t3.1}, Corollary~\ref{c3.1}, and Theorem~\ref{t3.2}
are closely related to a theorem of A.~D.~Aleksandrov about the
rigidity of the boundary
$\partial U$
of a strictly convex domain
$U$
in Euclidean
$n$-space
$\mathbb R^n$
by the relative metric
$\rho_{\partial U,U}$
on the boundary. The following is an important particular case of
this theorem:

\begin{theorem}\label{t3.3} {\rm(A.~D.~Aleksandrov
(see,~\cite{VA1})).}
Let
$U_1$
be a strictly convex domain in
$\mathbb R^n$
{\rm(}i.e.{\rm,} for any
$\alpha,\beta \in \cl U_1$ every shortest path
$\gamma = \gamma_{\alpha,\beta,\cl U_1} : [0,1] \to \cl U_1$
between
$\alpha$
and
$\beta$
{\rm(}in the metric
$\rho_{\cl U_1})$
satisfies
$\gamma(]0,1[) \subset U_1)$.
Assume that
$U_2 \subset \mathbb R^n$
is any domain whose closure is a Lipschitz manifold {\rm(}such
that
$\partial (\cl U_2) = \partial U_2 \ne\varnothing);$
moreover{\rm,}
$\partial U_1$
and
$\partial U_2$
are isometric in their relative metrics
$\rho_{\partial U_1,U_1}$
and
$\rho_{\partial U_2,U_2}$.
Then
$\partial U_1$
and
$\partial U_2$
are isometric in the Euclidean metric.
\end{theorem}

We say that an
$n$-dimensional
compact (closed) connected
$C^0$-submanifold
$Y$
with boundary
$\partial Y \ne\varnothing$
of an
$n$-dimensional
smooth connected
(respectively, $n$-dimensional
smooth complete connected) Riemannian
manifold
$X$
without boundary has property
${\rm(\circ)}$
if $\gamma_{x,y,Y}(]0,1[) \subset \inter Y$ for any two points
$x,y \in \partial Y$
and for every shortest path
$\gamma_{x,y,Y} : [0,1] \to Y$
in the metric
$\rho_Y$
joining these points.

\begin{theorem}\label{t3.4}
Let $n = 2$. Suppose that ${\rm(i)}$ holds for a $2$-dimensional
compact connected $C^0$-submanifold $Y_1$ with boundary $\partial
Y_1 \ne\varnothing$ in a $2$-dimensional smooth connected
Riemannian manifold $X$ without boundary{\rm;} moreover{\rm,}
$Y_1$ has property~${\rm(\circ)}$. Assume that $Y_2 \subset X$
is a $2$-dimensional compact connected $C^0$-submanifold
with $\partial Y_2 \ne\varnothing$
in
$X$
and $\partial
Y_1$ and $\partial Y_2$ are isometric in the metrics $\rho_{Y_j}$
$(j = 1,2)$. Then $\partial Y_2$ also has property~${\rm(\circ)}$.
\end{theorem}

This theorem has the following generalization.

\begin{theorem}\label{t3.5}
Let
$n = 2$.
Suppose that
${\rm(i)}$
holds for a
$2$-dimensional
closed connected
$C^0$-submanifold
$Y_1$
with boundary
$\partial Y_1\,\, (\ne\varnothing)$
in a
$2$-dimensional
smooth complete connected Riemannian manifold
$X$
without boundary satisfying
${\rm(\circ)}$.
Assume that
$Y_2 \subset X$
is a
$2$-dimensional
closed connected
$C^0$-submanifold
with
$\partial Y_2 \ne\varnothing$
in
$X;$
moreover{\rm,}
$\partial Y_1$
and
$\partial Y_2$
are isometric in the metrics
$\rho_{Y_j}$
$(j = 1,2)$.
Then
$Y_2$
has the property
${\rm(\circ)}$
as well.
\end{theorem}

\begin{corollary}\label{c3.2}
{\rm(of Theorem~\ref{t3.4}).}
Assume that the hypothesis of Theorem{\rm~\ref{t3.4}} hold and
that the manifold
$X$
has the following property{\rm:}
$\rho_X(x,y) = \rho_Y(x,y)$
for any two points
$x$
and
$y$
on the boundary
$\partial Y$
of every
$2$-dimensional
compact connected
$C^0$-submanifold
$Y \subset X$
with
$\partial Y \ne\varnothing$
satisfying
${\rm(i)}$
and
${\rm(\circ)}$.
Then
$\partial Y_1$
and
$\partial Y_2$
are isometric in the metric
$\rho_X$
of the ambient manifold~$X$.
\end{corollary}

\begin{corollary}\label{c3.3} {\rm(of Theorem~\ref{t3.5})}{\bf.}
Assume that the hypothesis of Theorem{\rm~\ref{t3.5}} hold and
that the manifold
$X$
has the following property{\rm:}
$\rho_X(x,y) = \rho_Y(x,y)$
for any two points
$x$
and
$y$
on the boundary
$\partial Y$
of every
$2$-dimensional
closed connected
$C^0$-submanifold
$Y \subset X$
with
$\partial Y \ne\varnothing$
satisfying
${\rm(i)}$
and
${\rm(\circ)}$.
Then
$\partial Y_1$
and
$\partial Y_2$
are isometric with respect to~$\rho_X$.
\end{corollary}

\begin{theorem}\label{t3.6}
Let
$n \ge 3$.
Suppose that
$(X,g)$
is an
$n$-dimensional
smooth connected Riemannian manifold whithout boundary and
$Y_1$
and
$Y_2$
are
$n$-dimensional
compact connected
$C^0$-submanifolds
with nonempty boundaries
$\partial Y_1$
and
$\partial Y_2$
in
$X$
satisfying conditions
${\rm(i)}$
and
${\rm(ii)}$
{\rm(}in Theorem{\rm~\ref{t3.2}}{\rm)}. Assume that
$Y_1$
has property~${\rm(\circ)}$
and~$X$ satisfies the following condition{\rm:}
$\rho_X(x,y) = \rho_Y(x,y)$
for any two points
$x$
and
$y$
on the boundary
$\partial Y$
of every
$n$-dimensional
compact connected
$C^0$-submanifold
$Y \subset X$
with
$\partial Y \ne\varnothing$
satisfying
${\rm(i)}$,
${\rm(ii)}$,
and
${\rm(\circ)}$.
Suppose also that
$\partial Y_1$
and
$\partial Y_2$
are isometric in the metrics
$\rho_{Y_j},$
where
$j = 1,2$.
Then
$\partial Y_1$
and
$\partial Y_2$
are isometric in~$\rho_X$.
\end{theorem}

\begin{theorem}\label{t3.7}
Let
$n \ge 3$.
Suppose that
$(X,g)$
is an
$n$-dimensional
smooth complete connected Riemannian manifold without boundary
and
$Y_1$
and
$Y_2$
are
$n$-dimensional
closed connected
$C^0$-submanifolds
with nonempty boundaries
$\partial Y_1$
and
$\partial Y_2$
in
$X$
satisfying
${\rm(i)}$
and
${\rm(ii)}$.
Assume that
$\partial Y_1$
has property~${\rm(\circ)}$
and~$X$
satisfies the following condition{\rm:}
$\rho_X(x,y) = \rho_Y(x,y)$
for any two points
$x$
and
$y$
on the boundary
$\partial Y$
of every
$n$-dimensional
closed connected
$C^0$-submanifold
$Y$
with
$\partial Y \ne\varnothing$
in
$X$
satisfying
${\rm(i)}$,
${\rm(ii)}$,
and
${\rm(\circ)}$.
Suppose also that
$\partial Y_1$
and
$\partial Y_2$
are isometric in the metrics
$\rho_{Y_j}$
$(j = 1,2)$.
Then
$\partial Y_1$
and
$\partial Y_2$
are isometric in~$\rho_X$.
\end{theorem}

\textit{Proofs of Theorems{\rm~\ref{t3.2}}
and{\rm~\ref{t3.4}}--{\rm\ref{t3.7}}}
are similar to the proof of Theorem~\ref{t3.1} (Theorems~\ref{t3.2}
and~\ref{t3.4}--\ref{t3.7} can be proved using the corresponding
analogs of Theorems~\ref{t2.1} and~\ref{t2.3}).
\vskip4mm

In conclusion, note that main results of our article were
earlier announced in~\cite{Ko} and~\cite{Ko1}.
\vskip3mm

{\bf Acknowledgements}
\vskip3mm

The authors were partially supported by the~Interdisciplinary
Project of the Siberian and Far-Eastern Divisions of the Russian
Academy of Sciences (2012-2014 no. 56), the State Maintenance
Program for the Leading Scientific Schools of the Russian
Federation (Grant NSh-921.2012.1), and the Exchange Program between
the Russian and Polish Academies of Sciences (Project 2014-2016).

\end{document}